\documentclass[12pt]{article}

\usepackage{amssymb}
\usepackage{amsmath,amssymb,amsthm, amscd}
\usepackage{indentfirst}

\theoremstyle{definition}

\title{\textbf{
   Regular Blaschke
  Para-Umbilical Hypersurfaces in the Conformal Space ${\mathbb Q}^n_s$}
}
 \author {{\bf Tongzhu Li$^1$, Changxiong Nie$^2$} \\
\small{1 Department of Mathematics, Beijing Institute of
Technology,} \\
\small{Beijing, China,100081, E-mail: litz@bit.edu.cn} \\
\small{2 Faculty of Mathematics and Statistics, Hubei Key Laboratory of Applied Mathematics,  }\\
\small{Hubei University, Wuhan, China, 430062, E-mail: chxnie@163.com
}}

\begin{document}

\maketitle \footnotetext [1]{T. Z. Li is   supported by the
grant No. 11571037  of NSFC; } \footnotetext [2]{C. X. Nie is
partially supported
by the grant No. 11571037  of NSFC.} 

\begin{abstract}
   In [2] we have classified the Blaschke
  quasi-umbilical submanifolds  in the conformal space ${\mathbb Q}^n_s$.
 In this paper we shall classify the Blaschke
  para-umbilical hypersurfaces  in the conformal space ${\mathbb Q}^n_s$.
 That may be also considered as the extension of the classification of the conformal
 isotropic submanifolds in the conformal space ${\mathbb Q}^n_s$.
\end{abstract}



\par\noindent
{\bf {\S} 1. Introduction.}
\par\medskip

Let $\mathbb R^{N}_{s}$ denote pseudo-Euclidean space, which is the real vector space $\mathbb R^{N}$
with the non-degenerate inner product $\langle,\rangle$ given by
  $$\langle\xi,\eta\rangle
  = - \sum_{i=1}^{ s } x_iy_i  +  \sum_{ i =  s + 1 }^{N} x_iy_i,
   $$
  where $\xi=(x_1, \cdots
x_{ _{N} } ),\eta=(y_1, \cdots , y_{ _{N} } )\in\mathbb R^{N}$.

Let
$$C^{n+1}:=\{\xi\in\mathbb R^{n+2}_{s+1}|\langle \xi  ,\xi  \rangle=0,\xi  \neq0\}, $$
 $$\mathbb Q^n_s:=\{[\xi  ]\in\mathbb R P^{n+1}|\langle \xi  ,\xi  \rangle=0\}=
 C^{n+1}/(\mathbb R\backslash\{0\}). $$
 We call
$C^{n+1}$ the light cone in $\mathbb R^{n+2}_{s+1}$ and  $\mathbb Q^n_s$ the
conformal space ({\sl or} projective light cone) in $\mathbb R P^{n+1}$.

  The standard metric $h$ of the conformal space $\mathbb Q^n_s$ can be obtained
 through the pseudo-Riemannian submersion
 $$\pi:C^{n+1}\rightarrow\mathbb Q^n_s,\xi\mapsto[\xi].$$
   We can check $(\mathbb Q^n_s,h)$ is a pseudo-Riemannian manifold.

  We define the pseudo-Riemannian sphere space $\mathbb S^n_s ( r ) $  and
pseudo-Riemannian hyperbolic space  $\mathbb H^n_s ( r ) $ with radius $r$ by
$$\mathbb S^n_s ( r ) =\{u\in \mathbb R^{n+1}_s|\langle u, u\rangle = r^2\},  \quad
\mathbb H^n_s ( r ) =\{u\in \mathbb R^{n+1}_{s+1}|\langle u, u\rangle = - r^2\}.   $$
When   $ r =1 $ we usually omit the radius $ r   $. When $s=1$ and $ r =1 $ we call them de Sitter space
$\mathbb S^n_1$ and anti-de Sitter space $\mathbb H^n_1$.

We may assume    ${\mathbb Q}^n_s$ to be the common compactification of
$\mathbb R^n_s$, $\mathbb S^n_s$ and $ \mathbb H^n_s$, and
 $\mathbb R^n_s$, $\mathbb S^n_s$ and $\mathbb H^n_s$ to be the subsets
 of $\mathbb
 Q^n_s$ when referring only to  the conformal geometry.

  When $s=0$, our analysis in this text can be reduced to the Moebius submanifold geometry in the
  sphere space (see  [4]).

This paper is organized as follows. In Section 2 we
 recall the  submanifold theory
 in the conformal space ${\mathbb Q}^n_s$ and give the relations  between conformal invariants and
 isometric
ones for submanifolds in several particular surroundings. In
Section 3  we  classify the conformal surfaces in   ${\mathbb Q}^3_1$. In
Section 4  we  classify    the Blaschke
   para-umbilical hypersurfaces in   ${\mathbb Q}^n_s$.\\

\bigskip

\par\noindent
{\bf {\S} 2. Fundamental equations.}
\par\medskip

 We recall the scheme of submanifold theory
 in the conformal space ${\mathbb Q}^n_s$ first.
 A classical theorem tells us that

 {\bf Theorem 2.1.}(see [3]) The conformal group of the conformal space
 ${\mathbb Q}^n_s$ is $O(n-s+1,s+1)/
 \{\pm 1\}$. If $\varphi$ is a conformal
 transformation on $\mathbb Q^n_s$,  then there is $A\in O(n-s+1,s+1)$,
 such that $\varphi=\Phi_A$ and $\Phi_A([X])=[XA]$.

Suppose that $x:\mathbf M^m_t\rightarrow{\mathbb Q}^n_s(s\geq1)$ is an $m$-dimensional
 Riemannian or pseudo-Riemannian submanifold
 with index $t(0\leq t\leq s)$. That
is, $x_*(\mathrm{T}\mathbf M)$ is non-degenerate subbundle of
$(\mathrm{T}{\mathbb Q}^n_s,h)$ with index $t(0\leq t\leq s)$.
  When $t=0$  we call $\mathbf M$ space-like submanifold. When $t>0$  we call $\mathbf M$
  pseudo-Riemmanian submanifold. Especially when $t=1$,
  $\mathbf M$ is called Lorentzian  submanifold or time-like submanifold.
  From now on, we always assume that the submanifold $x$  has
   index $t(0\leq t\leq s)$.

  Let $y:U\rightarrow C^{n+1}$ be a lift
of $x:\mathbf M\rightarrow{\mathbb Q}^n_s$ defined on an open subset
$U$ of $\mathbf M$. We denote by $\Delta$ and $\rho$ the Laplacian
operator and the scalar curvature of the local non-degenerate
metric $\langle\text dy, \text dy\rangle$. Then we have

{\bf Theorem 2.2.} (cf. [3]) Suppose that $x:\mathbf M\rightarrow{\mathbb Q}^n_s$ is an $m$-dimensional
 Riemannian or pseudo-Riemannian submanifold
  with index $t(0\leq t\leq s)$.
  On $\mathbf M$ the 2-form
  $$g:= \pm (\langle\Delta y,
\Delta y\rangle- \frac{ m }{ m - 1 }\rho)\langle\text dy, \text dy\rangle  $$ is a
  globally defined conformal invariant of $x $.

 {\bf Definition 2.1.} We call an
$m$-dimensional submanifold $x:\mathbf M\rightarrow {\mathbb Q}^n_s$
a regular submanifold if the 2-form $g:= \pm(\langle\Delta y, \Delta
 y\rangle- \frac{ m }{ m - 1 }\rho)\langle\text dy, \text dy\rangle$ is non-degenerate.
    $g$ is called the conformal metric of the regular
submanifold
 $x:\mathbf M\rightarrow {\mathbb Q}^n_s$.

 In this paper we assume that $x:\mathbf M\rightarrow {\mathbb Q}^n_s$
 is a regular  submanifold. Since the metric $g$ is
non-degenerate (we call it the conformal metric), there exists a
unique lift $Y:\mathbf M\rightarrow C^{n+1}$ such that
$g=\langle\text dY,\text dY\rangle$ up to sign. We call $Y$
the canonical lift of $x$.

{\bf Definition 2.2.} The two submanifolds $x,\tilde{x}$ are
conformally equivalent, if there exists a conformal transform
 $\sigma: {\mathbb Q}^n_s\rightarrow
{\mathbb Q}^n_s$, such that $\tilde{x}=\sigma\circ x$.

 It follows from   Theorem 2.1   that

{\bf Theorem 2.3.} Two submanifolds $x,\tilde{x}:\mathbf M\rightarrow
{\mathbb Q}^n_s$
 are
conformally equivalent if and only if there exists $T\in O(n-s+1,s+1)$ such
that $\tilde Y= T Y $, where $Y,\tilde Y$ are canonical lifts of
$x,\tilde x$
, respectively .\\

Let $\{e_1, \cdots , e_m\}$ be a local   basis of $\mathbf
M$ with dual basis $\{\omega^1, \cdots , $ $\omega^m\}$.  Denote
$Y_i=e_i(Y)$. We define
 $$N:=-\frac{1}{m}\Delta Y-\frac{1}{2m^2}\langle\Delta Y,
 \Delta Y\rangle Y.  $$
 Analogous to the corresponding calculation of  [13], we have
 $$\langle N, Y\rangle=1, \langle N, N\rangle=0,
 \langle N, Y_k\rangle=0,\quad1\leq k\leq m. $$
We may decompose $\mathbb R^{n+2}_{s+1}$ such that
$$\mathbb R^{n+2}_{s+1}=\text{span}\{Y, N\}\oplus \text{span}
\{Y_1, \cdots , Y_m\}\oplus\mathbb V  $$ where $\mathbb
V\bot\text{span}\{Y, N, Y_1, \cdots , Y_m\}$.  We call $\mathbb V$
 the conformal normal bundle for $x:\mathbf M\rightarrow {\mathbb Q}^n_s$.
 Let $\{\xi_{m+1},  \cdots , \xi_n\}$ be a local basis of the
bundle $\mathbb V$ over $\mathbf M$. Then $\{Y, N, Y_1, \cdots ,
 Y_m, \xi_{m+1},  \cdots , \xi_n\}$ forms a moving frame in
 $\mathbb R^{n+2}_{s+1}$ along $\mathbf M$. We adopt the conventions on the ranges
of indices in this paper without special claim:
  $$1\leq i, j, k, l,p,q\leq m;\quad m+1\leq\alpha, \beta,\gamma,\nu\leq n
  . $$

 We may write the
 structural equations as follows
 $$\mathrm{d}Y=\sum_i\omega^iY_i;\quad
 \mathrm{d}N=\sum_i\psi^iY_i+\sum_\alpha\phi^\alpha \xi_\alpha ; \eqno(2.1) $$
 $$\mathrm{d}Y_i=-\psi_iY-
 \omega_iN+\sum_j\omega_{i}^{j}Y_j+\sum_\alpha\omega^\alpha_i
 \xi_\alpha ; \eqno(2.2)$$
 $$\mathrm{d}\xi_\alpha =
 -\phi_\alpha Y     +  \sum_i\omega_{\alpha}^iY_i+\sum_\beta\omega^\beta_\alpha
 \xi_\alpha , \eqno(2.3) $$
 where the coefficients of $\{Y,N,Y_i,\xi_\alpha \}
$ are 1-forms on $\mathbf M$.

  It is clear that
 $\mathbb A:=\sum_i\psi_i\otimes\omega^i,
\mathbb
B:=\sum_{i,\alpha}\omega^\alpha_i\otimes\omega^i\xi_\alpha,\Phi:=\sum_\alpha\phi^\alpha
\xi_\alpha $ are globally defined conformal invariants. Let
$$
\psi_i=\sum_jA_{ij}\omega^j,
\quad\omega^\alpha_i  =   \sum_jB^\alpha_{ij}\omega^j,
\quad  \phi^\alpha  =   \sum_{i  }C^\alpha_i\omega^i.
 $$
  Denote the covariant derivatives of
these tensors  with respect to conformal metric
$g$ as follows:
$$\sum_jC^\alpha_{i, j}\omega^j=\mathrm dC^\alpha_i-\sum_jC^\alpha_j\omega_{i}^{j}+\sum_\beta
C^\beta_i\omega^\alpha_\beta;  $$
  $$\sum_kA_{ij, k}\omega^k =  \mathrm dA_{ij} - \sum_kA_{ik}\omega^k_j
  -  \sum_kA_{kj}\omega^k_i; $$
  $$\sum_kB^\alpha_{ij, k}\omega^k=\mathrm dB^\alpha_{ij}  -  \sum_kB^\alpha_{ik}\omega^k_j
  -  \sum_kB^\alpha_{kj}\omega^k_i+\sum_\beta
B^\beta_{ij}\omega^\alpha_\beta.  $$
 The curvature forms $\{ \Omega^{i}_{j} \}$ and the
 normal curvature forms $\{ \Omega^{\alpha}_{\beta} \}$ of the submanifold
 $x :\mathbf M\rightarrow
{\mathbb Q}^n_s$ can be written by
  $$ \Omega^{i}_{j}  =  \frac{1}{2}\sum_{kl}
  R^{i}_{\ jkl}\omega^k\wedge\omega^l
  =  \omega^{i}\wedge\psi_{j}  +  \psi^i\wedge\omega_j
  -   \sum_\alpha\omega^{i}_{\alpha}\wedge\omega^\alpha_j  ;
 $$
 $$\Omega^{\alpha}_{\beta}=  \frac{1}{2}\sum_{kl}
  R^{\alpha}_{\ \beta
  kl}\omega^k\wedge\omega^l
=  - \sum_i\omega^{\alpha}_{i}\wedge\omega^i_\beta
   . $$
     Denote
  $$ g_{ij}=\langle Y_i,Y_j\rangle,  \quad
 g_{\beta\gamma}=\langle \xi_\beta,
 \xi_\gamma\rangle, \quad
 (g^{ij})=(g_{ij})^{-1}, \quad
 (g^{\beta\gamma})=(g_{\beta\gamma})^{-1}, $$
 $$
 R _{ij
  kl}= \sum _p g_{ip}
  R^{p }_{\ j
  kl}, \quad
   R_{\alpha  \beta
  kl}= \sum_\nu g_{\alpha\nu}
  R^{\nu}_{\ \beta
  kl}.$$
 Then
the integrable conditions of the structure equations are
  $$A_{ij, k}-A_{ik, j}=
  -  \sum_{\alpha\beta}g_{\alpha\beta}
  (  B^\alpha_{ij}C^\beta_k  -  B^\alpha_{ik}C^\beta_j   ); \quad B^\alpha_{ij,
  k}-B^\alpha_{ik, j}=  g_{ij}C^\alpha_k - g_{ik}C^\alpha_j;
   $$
  $$C^\alpha_{i, j}-C^\alpha_{j, i}
  =\sum_{kl}g^{kl}(B^\alpha_{ik}A_{lj}-B^\alpha_{jk}A_{li});
  \quad R_{\alpha\beta
  ij}  =  \sum_{kl\gamma\nu}g_{\alpha\gamma}g_{\beta\nu}  g^{kl}
(B^\gamma_{ik}B^\nu_{lj}-B^\nu_{ik}B^\gamma_{lj});  $$
  $$R_{ijkl}=\sum_{\alpha\beta}g_{\alpha\beta}
(B^\alpha_{ik}B^\beta_{jl}-B^\alpha_{il}B^\beta_{jk})+
(g_{ik}A_{jl}-g_{il}A_{jk}) +(A_{ik}g_{jl}-
A_{il}g_{jk}). $$
 Furthermore, we have
  $$\text{tr}(\mathbb A)=\frac{1}{2m}(  \frac{ m }{ m - 1 }\rho\pm1);\quad
  R_{ij}=\text{tr}(\mathbb
  A) g_{ij}+(m-2)A_{ij}-\sum_{kl\alpha\beta}
  g^{kl}g _{\alpha\beta }B^\alpha_{ik}
  B^\beta_{lj} ; $$
  $$(1-m)C^\alpha_i=\sum_{jk}g^{jk}B^\alpha_{ij,k};\quad
  \sum_{ijkl\alpha\beta}g^{ij}g^{kl}
 g _{\alpha\beta}B^\alpha_{ik}B^\beta_{jl}=\pm\frac{m-1}{m};\quad
 \sum_{ij } g^{ij}B^\alpha_{ij}=0.  $$

 From above we know that in the case
$m\geq3$ all coefficients in the PDE system (2.1)-(2.3) are
determined by the conformal metric $g$,  the conformal second
fundamental form $\mathbb B$ and the normal connection
$\{\omega^\beta_\alpha\}$ in the conformal normal bundle. Then we have

{\bf Theorem 2.4.} Two hypersurfaces $x:\mathbf M^m_t\rightarrow
\mathbb Q^{m+1}_s$ and $\widetilde{x}:\widetilde{\mathbf
M}^m_t\rightarrow \mathbb Q^{m+1}_s(m\geq3)$ are conformal equivalent
if and only if there exists a diffeomorphism $f:\mathbf
M\rightarrow\widetilde{ \mathbf M}$ which preserves the conformal
metric   and the conformal second fundamental form. In another word,
$\{g, \mathbb B\}$ is a complete invariants system of the
hypersurface $x:\mathbf M^m\rightarrow
\mathbb Q^{m+1}_s(m\geq3)$.\\

 When $ \epsilon = 1, 0, - 1$, let the
 pseudo-Riemannian     space form $\mathbf R^n_s ( \epsilon ) $ denote
$ \mathbb S^n_s, \mathbb R^n_s$, $\mathbb H^n_s$, respectively.
 Let $\sigma_\epsilon:\mathbf R^n_s ( \epsilon )\rightarrow\mathbb{Q} ^n_s $
 be the standard conformal embedding( see [3]).

 Next we give the relations between the conformal
invariants induced above and isometric invariants of $u:\mathbf
M^m_t\rightarrow\mathbf R^n_s ( \epsilon )$.
 Let
 $\{e_1, \cdots , e_m\}$ be an local basis for $u$ with dual basis
 $\{\omega^1, \cdots , \omega^m\}$. Let $\{e_{m+1},  \cdots ,$ $ e_n\}$
 be a local basis of the normal bundle of $u$.
 Then we have the first and second fundamental forms $I, II$ and the
 mean curvature vector $\overrightarrow{ H }$.  We may write
 $$I  =  \sum_{ij}I_{ij}\omega^i\otimes\omega^j,\quad
 II=\sum_{ij\alpha}h^\alpha_{ij}\omega^i\otimes\omega^je_\alpha$$
 $$(I^{ij}) = (I_{ij})^{-1}, \quad
 \overrightarrow{ H } = \frac{1}{m}\sum_{ij\alpha}I^{ij}
 h^\alpha_{ij}e_\alpha:=\sum_\alpha H^\alpha e_\alpha.$$
  From the structure equations
 $$\mathrm{d}u=\sum_i\omega^iu_i,
  \quad\mathrm{d}u_i=\sum_j\theta^j_iu_j+\sum_\alpha\theta^\alpha_ie_\alpha
  -\epsilon\omega_i u,
  \quad\mathrm{d} e_\alpha
  =\sum_j\theta^j_\alpha u_j+\sum_\beta\theta^\beta_\alpha e_\beta,
 $$
  we have
 $$\Delta_I  u=m ( \overrightarrow{ H } - \epsilon u ),
 \quad\rho_I = m( m - 1 ) \epsilon +(m^2|\overrightarrow{ H }|^2
 -|II|^2), $$
 where
   $$|\overrightarrow{ H }|^2=\sum_{\alpha\beta}I_{\alpha\beta}H^\alpha H^\beta,
   I_{\alpha\beta}=(e_\alpha,e_\beta);\quad
    |II|^2=\sum_{ijkl\alpha\beta}
   I_{\alpha\beta}  I^{ik}  I^{jl}h^\alpha_{ij}h^{\beta}_{kl}.$$

 For the global lift
 $y:\mathbf M\rightarrow C^{n+1},$
    the conformal factor of $y$ is
  $$e^{2\tau} = \pm \frac{m}{m-1}(|II|^2-m|\overrightarrow{ H }|^2)
    .  \eqno{(2.4)}
  $$

     Furthermore, we have
      $$ \Delta_I  u=m(\overrightarrow{ H }-\epsilon u),
 \quad\rho_I = m^2|\overrightarrow{ H }|^2
 -|II|^2, \eqno{(2.5)}$$
 $$  A_{ij} = \tau_i\tau_j +
  \sum_{\alpha}h^\alpha_{ij}H_\alpha -\tau_{i,
  j}  -  \frac{1}{2}(\sum_{ij}I^{ij}\tau_i\tau_j +|\overrightarrow{ H }|^2 -\epsilon ) I_{ij}, \eqno{(2.6)}  $$
  $$  B^\alpha_{ij}  = e^{ \tau} ( h^\alpha_{ij}-H^\alpha I_{ij} ) ,
  \quad e^{ \tau} C^\alpha_i=  H^\alpha\tau_i
 -\sum_{j}h^\alpha_{ij}\tau^j-H^\alpha_{, i}
 ,  \eqno{(2.7)}$$
 where $\tau_{i, j}$ is the
Hessian of $\tau$ respect to $I$ and $H^\alpha_{, i}$ is the covariant
derivative of the mean curvature vector field of $u$ in the normal
bundle $N(\mathbf M)$ respect to $I$.

  \bigskip

\par\noindent
{\bf {\S} 3. Conformal surfaces in   ${\mathbb Q}^3_1$.}
\par\medskip

In this section let $x:\mathbf
  M^m_t\rightarrow {\mathbb Q}^{m+1}_s$ be an m-dimensional regular
  hypersurface with index
  $t(0\leq t\leq s)$. We use the notations in Section 2 and omit all normal
scripts in the formulas because the codimension now is one. Let
 $$ A^i_j = \sum_{k} g^{ik}A_{kj},\quad A=(A^i_j),$$
 $$B^i_j = \sum_{k} g^{ik}B_{kj},\quad B=(B^i_j).$$
 We rewrite some equations occurred preciously in the new form as follows
 $$\sum_{ij} B^{i}_{j}B_{i}^{j}= \frac{m-1}{m},\quad \sum_{i } B^{i}_{ i}=0,\eqno(3.1)$$
 $$  B_{i j,k}-B_{ik,j}= g_{ij}C_k-g_{ik}C_j,\quad
  A_{i j,k}-A_{ik,j}= B_{ij}C_k-B_{ik}C_j,\eqno(3.2)$$
  $$  C_{i, j }-C_{ j,i}=\sum_k( B_{i k}A^k_j-B_{j k}A^k_i),\eqno(3.3)$$
  $$\sum_{i } A^{i}_{ i}=\frac{1}{2m}(  \frac{ m }{ m - 1 }\rho\pm1).\eqno(3.4)$$

{\bf Definition 3.1.} We call an m-dimensional regular
submanifold $x:\mathbf
M\rightarrow {\mathbb Q}^n_s$  conformal if the conformal form
$   \Phi\equiv0. $

  Let $x:\mathbf
  M \rightarrow {\mathbb Q}^3_1$ be a regular space-like surface. We
  can write the structural equations as
  $$ e_i (N)=\sum_jA_i^jY_j + C_i\xi,\quad e_i(\xi)= C_i Y +\sum_jB_i^jY_j,\eqno(3.5)$$
 $$ e_j(Y_i)= - A_{ij}Y -g_{ij} N + \sum_k \Gamma^k_{ij} Y_k + B_{ij} \xi. $$

 Since
  $m=
2 $, we can find an
orthonormal basis $e_1,e_2$ of $x$ from (3.1) such that
  $$B=\mathrm{diag}(\frac{1}{2},-\frac{1}{2}).$$
  If $x$ is a conformal surface, we have $C_i=0,  i=1,2$.
  It implies from (3.2) that $B_{ij,k},A_{ij,k}$ are all symmetric
  with respect to the subscripts. For the same reason that $x$ has vanishing
  conformal form, by (3.3), we can modify the orthonormal basis $e_1,e_2$ such that
  $$A=\mathrm{diag}(a,b).$$
  Taking $i,j$  various values in
   $$\sum_kB _{ij, k}\omega^k=\mathrm dB _{ij}  -  \sum_kB _{ik}\omega^k_j
  -  \sum_kB _{kj}\omega^k_i, \eqno{( 3.6)}$$
  we have
  $$B_{11,i}=B_{22,i}=0,i=1,2.$$
  Therefore $B_{12,i}=0,i=1,2.$
  Letting $i=1,j=2$ in (3.6), we get the connection of $x$  is flat,
  {\it i.e.},
  $\omega^1_2=0$. It follows from (3.4) that
  $$ a+b=-\frac{1}{4} .\eqno(3.7)$$
  In addition, we may assume that there exist local  co-ordinates $u,v$ such that
  $$e_1=\frac{\partial}{\partial u},\quad
  e_2=\frac{\partial}{\partial v}.$$
Taking $i=1,j=2$   in
   $$\sum_kA _{ij, k}\omega^k=\mathrm dA _{ij}  -  \sum_kA _{ik}\omega^k_j
  -  \sum_kA _{kj}\omega^k_i, \eqno{( 3.8)}$$
 and noting $\omega^1_2=0$, we have
  $$A_{1 2,i}=0,i=1,2.$$
  So, when taking $i=j=1$ and $i=j=2$ in (3.8) respectively,
  we know that
  $a_v=b_u=0$.
  Adding (3.7) we shall see that $a,b$ are both constant.

  Next, we have the  structural equations as the following new form
  $$  N_u=aY_u,\quad N_v=bY_v,\quad
  \xi_u=\frac{1}{2}Y_u,\quad \xi_v=-\frac{1}{2}Y_v,\eqno(3.9)$$
 $$ Y_{uv}=0,\quad Y_{uu}=-aY-N+\frac{1}{2}\xi
 ,\quad Y_{vv}=-bY-N-\frac{1}{2}\xi.\eqno(3.10)$$
 So, we know from $Y_{uv}=0$ that Y can be split  as
 $$Y=F(u)+G(v).$$
 Substituting it into the  structural equations,
 we have
 $$F'''+(2a-\frac{1}{4})F'=0,\quad
  G'''+(2b-\frac{1}{4})G'=0.\eqno(3.11)$$
  By (3.7) we have
  $$(2a-\frac{1}{4}) +(2b-\frac{1}{4}) = -1.\eqno(3.12)$$
  In the following we discuss the resolve into three essential cases by noting the
 character of the  coefficients of the above PDEs (3.11).

  Case I: $2a-\frac{1}{4}<0,2b-\frac{1}{4}<0$.

  Let $2a-\frac{1}{4}=-r^2$. Then $2b-\frac{1}{4}=r^2-1$ and $0<r<1$. We have a particular resolve
  $$F=(r\cosh(ru), 0 , r\sinh (ru), 0 , 1 ),$$
  $$ G=(0,\sqrt{1-r^2}\cosh(\sqrt{1-r^2}v),0,
  \sqrt{1-r^2}\sinh(\sqrt{1-r^2}v),0).$$
  And we know that any resolve $(Y,N,Y_u,Y_v,\xi)$ of PDEs (3.9) and (3.10) is
  different from the initial resolve $(Y,N,Y_u,Y_v,\xi)_0$ up to a isometric transformation $T$ in $
  \mathbb R^5_2$, {\it i.e.}, $(Y,N,Y_u,Y_v,\xi)=T (Y,N,Y_u,Y_v,\xi)_0$.
  So $$Y=F+G:=(x,1)
   $$
   $$=(r\cosh(ru), \sqrt{1-r^2}\cosh(\sqrt{1-r^2}v)
   , r\sinh (ru), \sqrt{1-r^2}\sinh(\sqrt{1-r^2}v) , 1 )$$
   locally determines a surface $x:\mathbb H^1(r)\times\mathbb H^1(\sqrt{1-r^2})
   \rightarrow\mathbb H^3_1$
  whose canonical lift is $Y$.

 Case II: $2a-\frac{1}{4}<0,2b-\frac{1}{4}>0$.

  Let $2a-\frac{1}{4}=-r^2-1$. Then $2b-\frac{1}{4}=r^2 $ and $ r>0$. We have a particular resolve
  $$F=(1,r\cosh(\sqrt{r^2+1}u), r\sinh (\sqrt{r^2+1}u), 0 , 0  ),$$
  $$ G=(0,0,0,\sqrt{r^2+1}\cos(rv) ,
 \sqrt{r^2+1}\sin(rv) ).$$
  And we know that any resolve $(Y,N,Y_u,Y_v,\xi)$ of PDEs (3.9) and (3.10) is
  different from the initial resolve $(Y,N,Y_u,Y_v,\xi)_0$ up to a isometric transformation $T$ in $
  \mathbb R^5_2$, {\it i.e.}, $(Y,N,Y_u,Y_v,\xi)=T (Y,N,Y_u,Y_v,\xi)_0$.
  So $$Y=F+G:=(1,x)
   $$
   $$=(1,r\cosh(\sqrt{r^2+1}u), r\sinh (\sqrt{r^2+1}u), \sqrt{r^2+1}\cos(rv) ,
 \sqrt{r^2+1}\sin(rv) )$$
   locally determines a surface $x:\mathbb H^1(r)\times\mathbb S^1(\sqrt{r^2+1})
   \rightarrow\mathbb S^3_1$
  whose canonical lift is $Y$.

 Case III: $2a-\frac{1}{4}=-1,2b-\frac{1}{4}=0$.

  We have a particular resolve
  $$F=(0,\cosh u,\sinh u, 0,0),$$
  $$ G=(\frac{v^2}{2},0,0,   v ,
   \frac{v^2}{2}-1  ).$$
  And we know that any resolve $(Y,N,Y_u,Y_v,\xi)$ of PDEs (3.9) and (3.10) is
  different from the initial resolve $(Y,N,Y_u,Y_v,\xi)_0$ up to a isometric transformation $T$ in $
  \mathbb R^5_2$, {\it i.e.}, $(Y,N,Y_u,Y_v,\xi)=T (Y,N,Y_u,Y_v,\xi)_0$.
  So $$Y=F+G:=(\frac{\langle x,x\rangle}{2},x ,
   \frac{\langle x,x\rangle}{2}-1 )
   $$
   $$=(\frac{v^2}{2},\cosh u,\sinh u,   v ,
   \frac{v^2}{2}-1 )$$
   locally determines a surface $x:\mathbb H^1  \times\mathbb R^1
   \rightarrow\mathbb R^3_1$
  whose canonical lift is $Y$.

Summing up, we obtain

{\bf Theorem 3.1} If $x:\mathbf M^2\rightarrow\mathbb Q^3_1$
 is a space-like conformal surface, then it must be locally
 conformally equivalent to one of the three standard embedding surfaces:
$ \mathbb H^1(r)\times\mathbb H^1(\sqrt{1-r^2})
   \subset\mathbb H^3_1$,
 $ \mathbb H^1(r)\times\mathbb S^1(\sqrt{r^2+1})
   \subset\mathbb S^3_1$, and $  \mathbb H^1  \times\mathbb R^1 \subset\mathbb R^3_1$,
   where all radii of sphere or hyperbolic forms should be positive.

  Similarly, we shall get

{\bf Theorem 3.2} If $x:\mathbf M_1^2\rightarrow\mathbb Q^3_1$
 is a time-like conformal surface, then it must be locally
 conformally equivalent to one of the five standard embedding surfaces:
$ \mathbb H^1_1(r)\times\mathbb H^1(\sqrt{1-r^2})
   \subset\mathbb H^3_1$, $ \mathbb S^1_1(r)\times\mathbb H^1(\sqrt{1+r^2})
   \subset\mathbb H^3_1$,
 $ \mathbb S^1_1(r)\times\mathbb S^1(\sqrt{1-r^2 })
   \subset\mathbb S^3_1$, $  \mathbb R_1^1  \times\mathbb S^1 \subset\mathbb R^3_1$,
    and $  \mathbb S_1^1  \times\mathbb R^1 \subset\mathbb R^3_1$,
   where all the radii of (pseudo-Remannian) sphere  or hyperbolic forms should be positive.

  \bigskip

\par\noindent
{\bf {\S} 4.   Blaschke
  para-umbilical hypersurfaces in ${\mathbb Q}^n_s$.}
\par\medskip

 We remind readers that we shall retain the assumption on the head of Section 2.
 First, we give the

{\bf Definition 4.1.} We call an m-dimensional regular
hypersurface $x:\mathbf
M^m_t\rightarrow {\mathbb Q}^{m+1}_s$  Blaschke para-umbilical if there exist
a smooth function $\lambda,\mu$   on $\mathbf M$ such that
$$  A = \lambda I_m  + \mu
 B,  \ \text{ and}\ \  \Phi\equiv0,\eqno{(4.1)}$$
  where $I_m$ means $m$ order unit matrix.

{\bf Remark 4.1.}
  This definition is well-defined and it has no matter with the
  choose of local basis of $\mathbf M$.  When $ n = m + 1 $, a  Blaschke quasi-umbilical
submanifold reduces to a Blaschke para-umbilical hypersurface(c.f. [2]).

We have

{\bf Proposition 4.1.} If $u:\mathbf M^m_t\rightarrow \mathbf R^{m+1}_s ( \epsilon ) $ is a
  regular hypersurface with constant scalar curvature
  $ \rho_I $ and  mean curvature   $ H $,
then $x=\sigma_\epsilon \circ u$ is a  Blaschke para-umbilical hypersurface in
${\mathbb Q}^{m+1}_s$.

{ \sl Proof } Because of  (2.4) and (2.5),
 we know immediately that $ |  \overrightarrow{ H } | ^2$
 and $ | II |^2 $ are both constant.
 And one can easily see that the conformal factor $ e^{ 2 \tau }
 =\pm \frac{m}{m-1}(|II|^2-m|\overrightarrow{ H }|^2) =$constant.
 If the unit normal vector of hypersurface $x$ is space-like (or time-like), then
 we denote $\varepsilon=1$ (or $-1$). By use of (2.6)and (2.7), it follows from above that
$$  e^{2\tau}  A  =
   \varepsilon Hh + \frac{1}{2}(\epsilon-\varepsilon H^2) I_m, $$
  $$ e^{ \tau} B   =  h -H I_m   ,\quad
  C _i= 0,\forall i
. $$

If we choose
 $ \lambda = \frac{ 1 }{ 2 }e^{ -2 \tau }(   \epsilon +
 \varepsilon H^2 )$, and $\mu=\varepsilon e^{\tau} H$,
 we can verify that all the conditions of a Blaschke quasi-umbilical submanifold
are satisfied.$\Box$

{\bf Proposition 4.2.} Suppose that
$x:\mathbf
M^m_t\rightarrow {\mathbb Q}^{m+1}_s$   is a  Blaschke para-umbilical hypersurface in
${\mathbb Q}^{m+1}_s$.
Then the smooth function $ \lambda   $ in (4.1) must be constant.

{ \sl Proof } Suppose that $ \xi  $ is the unit normal vector of hypersurface $x$.
 Then from (3.1) and (3.5) we get
$$ e_i(N)  = \lambda e_i(Y) + \mu e_i( \xi) .   $$
  That means,
 $$
 \text dN+\lambda \text dY+\mu \text d\xi=0,\eqno(4.2)
  $$
 which implies that
 $$
 \text d\lambda\wedge \text dY+\text d\mu\wedge \text d\xi=0
 .$$
 Letting $\lambda_i=e_i(\lambda),
 \mu_i=e_i(\mu)$, combining with (3.5)  and the vanishing conformal form,
 we have
 $$\sum_{ijk} \lambda_i \omega^i\wedge\omega^j\delta^k_jY_k
 +\sum_{ijk} \mu_i \omega^i\wedge\omega^jB^k_jY_k=0.$$
 Because of the linear independence of
 $\{Y_1, \cdots , Y_m\}$ and the Cartan's lemma,
  we have
 $$
 \lambda_i\delta_{ j}^k
 +\mu_iB_{ j}^k=\lambda_j\delta_{i}^k
 +\mu_jB_{ i}^k
 .\eqno{(4.3)}$$
 Because $x$ has vanishing
  conformal form, by (3.3), we can
  choose an appropriate orthonormal basis $\{e_1,\cdots,e_m\}$ such that
  $$A=\mathrm{diag}(a_i),\quad
  B =\text{diag}(b_i) .$$
 For (4.2), fixing $i$, letting $j=k$, and taking summation over $j$,
 it follows from (3.4) that
 $$
 \lambda_i-\frac{1}{m-1}\mu_ib_i=0
 .\eqno{(4.4)}$$
 Taking $i\neq j= k$ in (4.2), we get
 $$
 \lambda_i+\mu_ib_j=0,\ \ i\neq j
 .\eqno{(4.5)}$$
 From (4.4) and (4.5) we have
 $$
 \mu_i(b_j+\frac{1}{m-1}b_i)=0,\ \  i\neq j
 .\eqno{(4.6)}$$
  If $\mu_i$'s are all zero,
 it follows from (4.4) that $\lambda_i$'s are all zero. Then
 $\lambda,\mu$ are both constant over $\mathbf M$.

 On the contrary, if $\mu_i$'s are not all zero,
 without the loss of generality,
 we may assume that $\mu_1\neq0$, then combining (3.1) and
 $$
 b_i=-\frac{1}{m-1}b_1,$$
 we can adjust the orient of the unit normal vector $\xi$ such that
 $$
 b_1=\frac{m-1}{m},\ \ b_2=\cdots=b_m=-\frac{1}{m}.\eqno{(4.7)}$$
 In the following we adopt the conventions on the ranges
of indices
 $$
 2\leq\alpha,\beta\leq m.$$
  Taking $i,j$  various values in
   (3.6),
  we have
  $$B_{11,i}=B_{\alpha\beta,i}=0,\forall i.$$
  Therefore $B_{1\alpha,i}=0,\forall i.$
  Taking $i=1,j=\alpha$    in
   (3.6),
  we have
 $
 \omega_{1\alpha}=0
 .$
 Similarly as precious induction in Section 3, we have
 $$
 R_{1\alpha1\alpha}=0=
 \varepsilon b_1b_2+a_1+a_\alpha
 ,\eqno{(4.8)}$$
 where $\varepsilon=\langle\xi,\xi\rangle$.
 So we know that $A=(a_1)\oplus(a_2I_{m-1})$.
 By (4.1) we get
 $$
 a_1+a_2    =     2\lambda+(b_1+b_2)\mu
 .\eqno{(4.9)}$$
 Combining (4.7)-(4.9), we get
 $$
 2\lambda+\frac{m-2}{m}\mu=\varepsilon
 \frac{m-1}{m^2},$$
 Therefore
 $$
 2\lambda_1+\frac{m-2}{m}\mu_1=0.\eqno{(4.10)}$$
 Taking $i=1,j=2$ in (4.5), we have
 $$
 \lambda_1=\frac{1}{m}\mu_1.\eqno{(4.11)}$$
 Substituting (4.11) into (4.10), we get
 $$
 \mu_1=0.$$
 This is a contraction to the assumption $\mu_1\neq0$.
 So,
 if $\mathbf M$ is connected, then
 $\lambda=\text{constant}, \mu=\text{constant}.\Box$

 If we take trace of the first equation of (4.1),
 we will find by (3.4) that
 $$m\lambda=   \text{tr}(  A)=
 \frac{1}{2m }(  \frac{ m }{ m - 1 }\rho\pm1)
 = \frac{1}{2(m   - 1) }\rho\pm\frac{1}{2m }.
 \eqno{(4.12)}$$
 which implies that the conformal scalar curvature
$$\rho=\text{constant}.$$
   Using the structural equations in Section 2, we have
 $$-m N = \Delta Y + \mathrm{tr} (A)Y.  $$
 From (4.12), we get
  $$-mN  =    \Delta Y  +m \lambda Y . \eqno{(4.13)}$$
Therefore by Proposition 4.1 and (4.2) we can find a constant vector $\overrightarrow c\in \mathbb R^{n+2}_{s+1}$ such that
$$N  =  \lambda Y +\mu \xi +   \overrightarrow c. \eqno{(4.14)}$$

It follows from (4.13) and (4.14) that
$$\langle  \overrightarrow c,Y\rangle   =  1,
 \quad\langle  \overrightarrow c,\xi\rangle   =  -\varepsilon\mu^2,
 \quad \langle\overrightarrow c, \overrightarrow c\rangle=
 -2\lambda  +  \varepsilon\mu^2,  \eqno{(4.15)}$$
 where $\varepsilon=\langle\xi,\xi\rangle$.

Then we discuss into the following three cases.

 {\bf Case 1}: $\langle\overrightarrow c, \overrightarrow c\rangle  =
 -2\lambda +  \varepsilon\mu^2 =0$.

 By use of an isometric transform of  $\mathbb R^{n+2}_{s+1}$ if necessary,
 assume that
$$\overrightarrow c=( 1,\mathbf0,  1).  $$
Letting
$$Y=(x_1,u,x_{n+2}),  $$
 it follows from the first equation of (4.15) and $\langle Y,Y\rangle=0$ that
$$Y=(\frac{  \langle u , u  \rangle  -1}{2}, u, \frac{  \langle u , u  \rangle  +1}{2}).
$$
 Then $x$ determines a hypersurface $u:\mathbf M^m_t\rightarrow\mathbb
R^{m+1}_s$ with
  $$I=  \langle \text du , \text du  \rangle  =\langle\text dY, \text
dY\rangle=g,$$
    which implies that
$$\Delta_{I} = \Delta,\ \ \ \rho_{I}=\rho=\text{constant}.   $$
 We know from [2] that
$$\xi= \varepsilon H Y + (0,\zeta,0),\eqno(4.16)
$$
where $\zeta$ is the unit normal vector of $u$.
  It  follows from the first and the second equations of (4.15) and (4.16) that
  $$H=-\mu^2=\text{constant}.$$
  Then,  $u$ is a
regular hypersurface with
constant scalar
curvature and  mean curvature   in $\mathbb R^{m+1}_s$.  In this case $x$ is locally
conformally equivalent to   a regular hypersurface with
constant scalar
curvature and  mean curvature   in $\mathbb R^{m+1}_s$.

 {\bf Case 2}: $\langle\overrightarrow c, \overrightarrow c\rangle  =
 -2\lambda +  \varepsilon\mu^2  :=- r^2, r={\rm constant}>0. $

  By use of an isometric transform of  $\mathbb R^{n+2}_{s+1}$ if necessary,
 assume that
$$\overrightarrow c=(r,\mathbf{0}).  $$
 Letting
$$Y=( x_{1} , u/r),  $$
 by similar method as above  we have
$$x_{ 1 }=1/r.  $$
So
$$Y=(1,u)/r,  \quad \langle u , u  \rangle  = 1 .  $$
Then $x$ determines a hypersurface $u:\mathbf M^m_t\rightarrow\mathbb
S^{m+1}_s$ with
  $$I/r^2=  \langle \text du , \text du  \rangle  /r^2=\langle\text dY, \text
dY\rangle=g, $$
   which implies that
  $$  r^2  \Delta_{I} = \Delta,\ \ \ \rho_{I}=\rho/r^2=\text{constant}. $$
   We know from [2] that
$$\xi= \varepsilon H Y + (0,\zeta ),\eqno(4.17)
$$
where $\zeta$ is the unit normal vector of $u$.
  It  follows from the first and the second equations of (4.15) and (4.17) that
  $$H=-\mu^2=\text{constant}.$$
  Then,  $u$ is a
regular hypersurface with
constant scalar
curvature and  mean curvature   in $\mathbb S^{m+1}_s$.  In this case $x$ is locally
conformally equivalent to   a regular hypersurface with
constant scalar
curvature and  mean curvature   in $\mathbb S^{m+1}_s$.

 {\bf Case 3}: $\langle\overrightarrow c, \overrightarrow c\rangle  =
 -2\lambda  +  \varepsilon\mu^2 :=r^2,  r={\rm constant}>0. $

  By use of an isometric transform of  $\mathbb R^{n+2}_{s+1}$ if necessary,
 assume that
$$\overrightarrow c=( \mathbf{0},   r).  $$
Letting
$$Y=( u/r, x_{n+2}),  $$
similarly  we have
$$x_{n+2}=1/r.  $$
So
$$Y=( u,1)/r,   \langle u ,
  u  \rangle  =-1 .  $$
Then $x$ determines a hypersurface $u:\mathbf M^m_t\rightarrow\mathbb
H^{m+1}_s$ with
  $$I/r^2=  \langle \text du , \text du  \rangle  /r^2=\langle\text dY, \text
dY\rangle=g,$$
    which implies that
  $$  r^2  \Delta_{I} = \Delta,\ \ \ \rho_{I}=\rho/r^2=\text{constant}. $$
   We know from [2] that
$$\xi= \varepsilon H Y + ( \zeta,0),\eqno(4.18)
$$
where $\zeta$ is the unit normal vector of $u$.
  It  follows from the first and the second equations of (4.15) and (4.18) that
  $$H=-\mu^2=\text{constant}.$$
  Then,  $u$ is a
regular hypersurface with
constant scalar
curvature and  mean curvature   in $\mathbb H^{m+1}_s$.  In this case $x$ is locally
conformally equivalent to   a regular hypersurface with
constant scalar
curvature and  mean curvature   in $\mathbb H^{m+1}_s$.

 So combining    Proposition 4.1   we get

 {\bf
Theorem 4.1.} Any Blaschke para-umbilic  hypersurface in ${\mathbb Q}^n_s$ is locally conformally equivalent to a
 regular hypersurface with
constant scalar
curvature and  mean curvature   in $\mathbb R^n_s,\mathbb S^n_s$, or $\mathbb H^n_s$.

\bigskip

{\bf Acknowledgements}

The
authors is deeply grateful to  Professor Changping Wang's guidance and inspiration.

\end{document}